\documentclass[12pt]{amsart}
\usepackage[hmargin=1in,vmargin=1.25in]{geometry}

\usepackage[utf8]{inputenc}
\usepackage{amsmath}
\usepackage{amsthm}
\usepackage{amssymb}
\usepackage{amsfonts}
\usepackage{enumerate}
\usepackage{tikz}
\usepackage{tkz-euclide}
\usepackage{ytableau}
\usepackage{mathtools}

\theoremstyle{plain}
\newtheorem{theorem}{Theorem}[section]
\newtheorem{lemma}[theorem]{Lemma}
\newtheorem{corollary}[theorem]{Corollary}

\theoremstyle{definition}
\newtheorem{defin}[theorem]{Definition}
\newtheorem{remark}[theorem]{Remark}
\newtheorem{exmp}[theorem]{Example}

\def\bZ{{{\mathbb Z}}}

\def\bx{{{\mathbf x}}}
\def\bk{{{\mathbf k}}}

\def\bu{{{\mathbf {u}}}}

\ytableausetup{boxsize=1.5em}

\newcommand{\ZZ}{\ensuremath{\mathbb{Z}}}
\newcommand{\inc}{\text{\rm inc}}
\newcommand{\SSYT}{\text{\rm SSYT}}
\newcommand{\CT}{\text{\rm CT}}
\newcommand{\CSSYT}{\text{\rm CSSYT}}
\newcommand{\col}{\text{\rm col}}
\newcommand{\threeone}{\ensuremath{(\mathbf{3}+\mathbf{1})}}
\newcommand{\sgn}{\text{\rm sgn}}
\newcommand{\swap}{\text{\rm swap}}

\title{Cylindric  $P$-tableaux for 3+1-free posets}

   \author{Isaiah Siegl}

   \address
   {Dept.\ of Mathematics\\
    University of Washington\\
    Seattle, WA}
   \email{isaiahsiegl@gmail.com}

\begin{document}

\begin{abstract}
For a $(\mathbf{3}+\mathbf{1})$-free poset  $P$, we define a hybrid of  $P$-tableaux and cylindric tableaux called cylindric $P$-tableaux. We introduce $P$-analogs of cylindric Schur functions, defined by a determinantal formula, and prove that they are  the weight generating functions of cylindric $P$-tableaux.
We deduce that certain sums of the  $e$-expansion coefficients of the chromatic symmetric function  $X_{\text{\rm inc}(P)}$
are positive. This improves on Gasharov's theorem on the Schur positivity of $X_{\text{\rm inc}(P)}$ and gives further evidence for the Stanley-Stembridge conjecture.
\end{abstract}%

\maketitle

\section{Introduction}
The chromatic symmetric function  $X_G(\mathbf{x})$ of a graph $G$ is
the sum $\sum_{\kappa} \mathbf{x}_\kappa$
over all proper colorings $\kappa: V(G) \to \mathbb{Z}_{>0}$ where $\mathbf{x}_{\kappa} = \prod_{v \in V(G)} x_{\kappa(v)}$.
It is  of particular interest when $G$ is the incomparability graph  $\inc(P) $ of a \threeone-free poset  $P$. 
In this case, a theorem of Haiman showed that $X_{\inc(P)}(\mathbf{x})$ is Schur positive \cite{Himmanant}, and a combinatorial formula for the Schur expansion was given by Gasharov \cite{Gasharov}. Stanley and Stembridge \cite{Stanleychromatic, StanleyStembridge}
conjectured that $X_{\inc(P)}(\mathbf{x})$ is a  positive sum of elementary symmetric functions;
the coefficient of $e_\lambda$ is known to be positive when $\lambda$ is a rectangular \cite{stembridge, CHSS}, two-column \cite{Wolfgangthesis, CHSS, Hwang}, or hook shape \cite{Wolfgangthesis, Hwang}.
Additionally, $X_G$ is known to be $e$-positive for several classes of graphs \cite{gebhard2001chromatic, dahlberg2018lollipop, harada2019cohomology, cho2019positivity, Dahlberg, ChoHong}.

In the recent paper \cite{Hwang} and in a companion paper to this one \cite{BEPS},
the authors adapted the Fomin-Greene theory of noncommutative Schur functions \cite{FG} to study chromatic symmetric functions.
This is the starting point for our work here, so
we review the basic setup.
For a finite poset $P$, let  $\ZZ[\mathbf{u}]$ be the polynomial ring in the commuting variables
$\mathbf{u} = \{u_{p}: p \in P\}$.
We define the \emph{$P$-elementary function} $e_k^P(\mathbf{u}) \in \ZZ[\mathbf{u}]$ by
\begin{equation}
e_k^P(\mathbf{u}) = \sum_{i_1 <_P i_2 <_P \cdots <_P i_k} u_{i_1}u_{i_2}\cdots u_{i_k}.
\end{equation}
(\cite{Hwang} and \cite{BEPS} are largely concerned with the setting where the $u_p$ are noncommuting variables but we only need the
commuting setup here.)

Let  $\Lambda(\mathbf{x})$ denote the ring of symmetric functions in variables $\mathbf{x}= x_1,x_2,\dots$.
Now, for $f(\mathbf{x}) \in \Lambda(\mathbf{x})$, we define the $P$-analog $f^P(\mathbf{u})$ of $f(\mathbf{x})$ to
be  the image of  $f(\mathbf{x})$ under the homomorphism
\begin{equation}
\psi \colon \Lambda(\mathbf{x}) \to \ZZ[\mathbf{u}], \quad  e_k(\mathbf{x}) \mapsto e_k^P(\mathbf{u}).
\end{equation}
The following result is a consequence of \cite[Theorem 2.6 and Eq. (2.21)]{BEPS}.
\begin{theorem}
The  $e$-expansion of $X_{\inc(P)}(\mathbf{x})$ can be expressed in terms of the  $P$-analogs $m_\lambda^P(\mathbf{u})$
of the monomial symmetric functions $m_\lambda(\mathbf{x})$ as
\begin{equation}
\label{e X via m u}
X_{\inc(P)}(\mathbf{x}) = \sum_{\lambda} \langle {\mathbf{u}_P} \rangle m_{\lambda}^P(\mathbf{u}) \, e_{\lambda}(\mathbf{x}),
\end{equation}
where $\langle {\mathbf{u}_P} \rangle m_{\lambda}^P(\mathbf{u})$ is the coefficient of $\mathbf{u}_P = \prod_{p \in P} u_p$ in the  $\mathbf{u}$-monomial expansion of  $m_{\lambda}^P(\mathbf{u})$.


Hence the  $\mathbf{u}$-monomial positivity of  $m_{\lambda}^P(\mathbf{u})$ implies the
Stanley-Stembridge conjecture (in fact, they're equivalent by \cite[Remark 3.15]{BEPS}).
\end{theorem}

On the other hand, Gasharov's result \cite{Gasharov} can be rephrased as giving a  $\mathbf{u}$-monomial positive formula for  $s_\lambda^P(\mathbf{u})$ in terms of $P$-tableaux (see Theorem \ref{t Gasharov sP}).  Hence a natural intermediate goal to go from
the known $\mathbf{u}$-monomial positivity of $s_\lambda^P(\mathbf{u})$ to the conjectured
$\mathbf{u}$-monomial positivity of $m_\lambda^P(\mathbf{u})$ is to establish $\mathbf{u}$-monomial positivity
for various  $P$-analogs $\psi(f(\mathbf{x}))$ for symmetric  $f(\mathbf{x})$ which ``lie in between''  monomial symmetric functions and Schur functions.

One such class are the cylindric Schur functions  $s_{\lambda/\mu/d}(\mathbf{x})$, where ``lie in between'' has the precise meaning that they are positive sums of some subset of monomials which sum to a Schur function.
The cylindric Schur functions are based on cylindric partitions of Gessel and Krattenthaler \cite{GKcyl}, and were used by Postnikov to study the quantum cohomology of the Grassmannian \cite{PostnikovCylSchur}. Lam showed that cylindric Schur functions are special cases of skew affine Schur functions \cite{LamAffStan}. McNamara \cite{McNamara} conjectured and Lee \cite{LeeCylSchur} proved that cylindric skew Schur functions expand positively in terms of cylindric Schur functions and that the coefficients of this expansion are the same as 3-point Gromov-Witten invariants.

We define  $P$-analogs  $s^P_{\lambda/\mu/d}(\mathbf{u}) := \psi(s_{\lambda/\mu/d}(\mathbf{x}))$ of the cylindric Schur functions and prove that these also have a formula in terms of  \emph{cylindric $P$-tableaux}, a hybrid
of cylindric tableaux and $P$-tableaux (Theorem \ref{main theorem}).
In particular, this establishes that the $s^P_{\lambda/\mu/d}(\mathbf{u})$ are  $\mathbf{u}$-monomial positive and
that certain sums of the coefficients $c^P_\lambda$ in  $X = \sum_\lambda c^P_\lambda e_\lambda(\mathbf{x})$ are positive (Corollary \ref{c sum of e coefs}); these sums are \emph{smaller} than those obtained from Gasharov's result in the same fashion.
As a special case, we recover a positive formula of Stembridge \cite{StembridgeImmanant} and Clearman-Hyatt-Shelton-Skandera \cite{CHSS}
for the coefficient $c^P_\lambda$ when $\lambda$ is a rectangle.

\section{$P$-cylindric Schur functions are  $\mathbf{u}$-monomial positive}
\subsection{Preliminaries}
For an integer partition $\lambda = (\lambda_1 \ge \lambda_2 \ge \cdots\ge \lambda_l) \vdash n$, the (English style) Young diagram of shape $\lambda$ is the set $\{(i,j): 1 \le i \le l, 1 \le j \le \lambda_i\}$, drawn as boxes labeled with  matrix-style coordinates where the box $(i,j)$ is in row  $i$ and column $j$. We often identify partitions with their corresponding Young diagrams so that, for partitions $\lambda, \mu$, we write $\mu \subset \lambda$ to mean that the Young diagram of $\mu$ is contained in the Young diagram of $\lambda$.
For $\mu \subset \lambda$, the skew shape $\lambda/\mu$ is the difference of Young diagrams $\lambda - \mu$.
We write  $\lambda'$ for the transpose partition of  $\lambda$ and $\ell(\lambda)$ for the number of nonzero parts of $\lambda$.

A semistandard Young tableau of shape $\lambda/\mu$ is a function $T: \lambda/\mu \to \bZ_{> 0}$ so that $T(i,j) < T(i+1,j)$ and $T(i,j) \le T(i,j+1)$, i.e. a Young tableau is an assignment of positive integers to the boxes of  $\lambda/\mu$ so that the entries strictly increase down columns and weakly increase across rows. We write $\SSYT(\lambda/\mu)$ for the set of semistandard Young tableaux of shape $\lambda/\mu$.

\subsection{Reformulation of Gasharov's theorem}
To set the stage for our main theorem on  $P$-cylindric Schur functions, we discuss several
precursors of this result.

It is well know from
the theory of symmetric functions
that the following two quantities are equal; either can be taken as the definition of the
skew Schur function  $s_{\lambda/\mu}(\bx)$.
\begin{equation}
\sum_{T \in \SSYT(\lambda/\mu)} \bx^T = \det [e_{\lambda_i' - \mu_j' - i+j}(\bx)]_{i,j=1}^{\ell(\lambda')},
\end{equation}
where  $\bx^T = \prod_{b\in \lambda/\mu}x_{T(b)}$.
Thus, for a poset  $P$, the  $P$-analog of $s_{\lambda/\mu}(\mathbf{x})$ is given by
\begin{equation}
s_{\lambda/\mu}^P(\bu) := \psi(s_{\lambda/\mu}(\mathbf{x})) = \det [e_{\lambda_i' - \mu_j' - i+j}^P(\bu)]_{i,j=1}^{\ell(\lambda')}.
\end{equation}
We say $T: \lambda/\mu \to P$ is a \emph{$P$-tableau} of shape $\lambda/\mu$ if $T(i,j) <_P T(i+1,j)$ and $T(i,j) \not>_P T(i,j+1)$, i.e. an assignment of elements of $P$ to the boxes of $\lambda/\mu$ is a $P$-tableau if it is increasing in $P$ down columns and non-decreasing in $P$ across rows.
\begin{theorem}
\label{t Gasharov sP}
For a \threeone-free poset  $P$,
\begin{equation}
s_{\lambda/\mu}^P(\bu) = \sum_{T \in \SSYT_P(\lambda/\mu)} \mathbf{u}^T,
\end{equation}
where $\SSYT_P(\lambda/\mu)$ denotes the set of $P$-tableaux of shape $\lambda/\mu$ and
$\bu^T = \prod_{b\in \lambda/\mu}u_{T(b)}$.
\end{theorem}
For straight shapes ($\mu = \varnothing$), this is a reformulation of Gasharov's theorem \cite{Gasharov} (it also follows directly from \cite[Theorem 3.9]{BEPS}); for skew shapes, it will follow from the more general Theorem \ref{main theorem} below.

\subsection{$P$-cylindric Schur functions}

\begin{defin}[Cylindric $P$-tableaux]
For a poset $P$, a $P$-tableau $T$ of shape $\lambda/\mu$, and a nonnegative integer $d$, define $T^d$ to be the diagram obtained by gluing a copy of the first column of $T$ to the last column of $T$ so that the copied column is shifted up by $d$ cells. i.e. $T^d$ has shape $\nu/\theta$ where $\nu' = (\lambda_1' + d, \lambda_2' + d,\ldots,\lambda_l'+d, \lambda_1')$ and $\theta' = (\mu_1' + d, \mu_2' + d,\dots,\mu_l'+d, \mu_1')$.

A \emph{cylindric $P$-tableau} of shape $\lambda/\mu/d$ is a function  $T \colon \lambda/\mu \to P$ such that  $T^d$
is a $P$-tableau.
Let 
\begin{equation}
\CT_P(\lambda/\mu/d) = \ \text{the set of cylindric $P$-tableau of shape $\lambda/\mu/d$}.
\end{equation}
Note that $T^d$ will never be a $P$-tableau when $d < \lambda_1' - \lambda_{\lambda_1}' $ or $d <  \mu_1'-\mu'_{\lambda_1}$, and $\CT_P(\lambda/\mu/d) = \CT_P(\lambda/\mu/d')$ when $d,d' \ge \lambda_1'$. Therefore we take $\max(\lambda_1' - \lambda_{\lambda_1}', \mu_1'-\mu'_{\lambda_1}) \le d \le \lambda_1'$, and in this case we say  $\lambda/\mu/d$ is a \emph{cylindric shape}.
\end{defin}

\begin{exmp}
Let $P$ be the total order on $ \bZ_{>0}$ and $T$ be the following tableau
\begin{equation}
\begin{ytableau}
\none & \none & 1 & 2 & 3 \\
\none & 1 & 3 & 5 & 5 \\
2 & 2 & 4 & 6 & 6 \\
3 & 4 & 5 & 7
\end{ytableau}.
\end{equation}
Then
\begin{eqnarray}
\raisebox{-20pt}{$T^3 =$ }
\begin{ytableau}
\none & \none & \none & \none & \none & 2 \\
\none & \none & 1 & 2 & 3 & 3 \\
\none & 1 & 3 & 5 & 5 \\
2 & 2 & 4 & 6 & 6 \\
3 & 4 & 5 & 7
\end{ytableau} & \raisebox{-20pt}{ and } &
\raisebox{-20pt}{$T^2 =$ }\begin{ytableau}
\none & \none & 1 & 2 & 3 & 2\\
\none & 1 & 3 & 5 & 5 & 3\\
2 & 2 & 4 & 6 & 6 \\
3 & 4 & 5 & 7
\end{ytableau}.
\end{eqnarray}
\end{exmp}

When $P $ is the total order on $ \bZ_{>0}$, we call cylindric $P$-tableaux \emph{cylindric semistandard Young tableaux} and write $\CSSYT(\lambda/\mu/d) = \CT_P(\lambda/\mu/d)$.

For a cylindric shape $\lambda/\mu/d$ define the cylindric Schur function $s_{\lambda/\mu/d}(\mathbf{x})$ by
\begin{equation}
\label{e cyl schur}
s_{\lambda/\mu/d}(\mathbf{x}) = \sum_{T \in \CSSYT(\lambda/\mu/d)} \bx^T.
\end{equation}

\begin{remark}
Our definition of cylindric semistandard Young tableaux and cylindric Schur functions are related to other definitions in the literature as follows.

(i) Our definition of the shape of a cylindric tableau differs from the definition used in \cite{GKcyl}. Our definition is slightly more general. For a cylindric shape $\lambda/\mu/d$ in our notation with $d > 0$, we convert to a cylindric shape $\nu/\eta/m$ in the notation of Gessel and Krattenthaler by taking
\begin{equation}
\nu = (\lambda_1 + \lambda_{1+d} + \lambda_{1+2d}+\cdots, \lambda_{2} + \lambda_{2+d}+\lambda_{2+2d}+\cdots,\dots,\lambda_d + \lambda_{2d}+\cdots),
\end{equation}
\begin{eqnarray}
\eta = \mu, & \text{ and } & m = \lambda_1.
\end{eqnarray}

(ii) For  $d > 0$, our definition \eqref{e cyl schur} of $s_{\lambda/\mu/d}(\mathbf{x})$ is
equivalent to the definition used by Postnikov \cite{PostnikovCylSchur} and McNamara \cite{McNamara}, but our conventions for cylindric shapes are different. Following \cite[\S4]{McNamara}, we convert a cylindric shape $\lambda/\mu/d$ in our notation to a shape $\nu/m/\theta$ in the notation of Postnikov and McNamara in the following way:
let $k = d$ and $n = \lambda_1 + d$. Then $\theta = \mu$, and $\nu$ is obtained from $\lambda$ by removing $m$ $n$-ribbons from the border of $\lambda$, where $m$ is the smallest number so that $\nu$ has at most $k$ rows.

For example, for $\lambda = (4,4,4,4,2,1,1)$, $\mu = (2,1)$, and $d = 3$, to convert to the notation of Postnikov and McNamara, we repeatedly remove ribbons with 7 cells from the bottom edge of $\lambda$ until the resulting shape has at most $3$ rows. In the following figure, we label the cells removed in the first ribbon with $1$'s and label the cells in the second ribbon with $2$'s.
\begin{equation*}
{\tiny
\begin{ytableau}
\none & \none & & 2\\
\none & & & 2\\
2 & 2 & 2 & 2\\
2 & 1&1 & 1\\
1 & 1 \\
1\\
1\\
\end{ytableau}}
\end{equation*}
Therefore we have that the equivalent shape $\nu/m/\theta$ in the notation of Postnikov and McNamara is $\nu = (3,3)$, $m=2$, and $\theta = (2,1)$.
\end{remark}

Gessel-Krattenthaler \cite{GKcyl} describe a method for expressing cylindric Schur functions as a sum of determinants and give an explicit formula in the case when $\lambda$ is a rectangular shape and $d = 0$. Postnikov \cite[Eq. (11)]{PostnikovCylSchur} (see also \cite[\S6]{McNamara}) then gave the following Jacobi-Trudi-like identity for cylindric Schur functions, making the result of Gessel-Krattenthaler explicit for arbitrary $\lambda$ and $d > 0$:
\begin{theorem}
\label{t GK}
\begin{equation}
s_{\lambda/\mu/d}(\mathbf{x}) = \sum_{\substack{k_1+k_2+\cdots +k_{\lambda_1} = 0 \\k_i \in \bZ}} \det \begin{bmatrix} e_{k_{i}(\lambda_1 + d) + \lambda_i' - \mu_j' - i + j}(\mathbf{x}) \end{bmatrix}.
\end{equation}
\end{theorem}

Now define the \emph{$P$-cylindric Schur functions} by
\begin{equation}
s_{\lambda/\mu/d}^P(\mathbf{u}) := \psi(s_{\lambda/\mu/d}^P(\mathbf{x})) = \sum_{\substack{k_1+k_2+ \cdots +k_r = 0 \\k_i \in \bZ}} \det \begin{bmatrix} e_{k_i(\lambda_1 + d) + \lambda_i' - \mu_j' - i + j}^P(\mathbf{u}) \end{bmatrix}.
\end{equation}
Our main result is the following  $P$-version of
Theorem \ref{t GK}.
\begin{theorem}
\label{main theorem}
For any \threeone-free poset  $P$ and cylindric shape $\lambda/\mu/d$,
\begin{equation}
s_{\lambda/\mu/d}^P(\bu) = \sum_{T \in \CT_P(\lambda/\mu/d)} \bu^T.
\end{equation}
\end{theorem}

By \eqref{e X via m u}, this $\mathbf{u}$-monomial positive formula for the  $P$-cylindric Schur functions has the following consequence for the expansion of  $X_{\inc(P)}$ into elementary symmetric functions.
\begin{corollary}
\label{c sum of e coefs}
Letting  $a_\nu$ denote the coefficients in
the monomial expansion
$s_{\lambda/\mu/d}(\mathbf{x}) = \sum_\nu a_\nu m_\nu(\mathbf{x})$ and  $c^P_\lambda$ the coefficients
in  $X_{\inc(P)} = \sum_\lambda c^P_\lambda e_\lambda(\mathbf{x})$, we have
\begin{align}
\label{e c sum of e coefs}
\sum_\nu a_\nu c^P_\nu = \langle \mathbf{u}_P \rangle \,  s_{\lambda/\mu/d}^P(\bu) \, = \,
\#\big\{T \in \CT_P(\lambda/\mu/d) : \text{ $T$ is standard}\big\},
\end{align}
where  a cylindric $P$-tableau is \emph{standard} if it contains each element of  $P$ exactly once.
\end{corollary}

Note that in the special case  $s_{\lambda/\mu/d}(\mathbf{x}) = s_\lambda(\mathbf{x})$ (when  $\mu=\varnothing$, $d =\lambda'_1$),  $a_\nu$ is the Kostka coefficient  $K_{\lambda \nu}$ and \eqref{e c sum of e coefs} follows
directly from Gasharov's result \cite{Gasharov}.
Thus  \eqref{e c sum of e coefs} can be regarded as a strengthening of this result as the  $a_\nu$ are typically smaller than Kostka coefficients.

See \S\ref{s special cases} for two other interesting cases of \eqref{e c sum of e coefs}.

\section{Proof of Theorem \ref{main theorem}}

Our proof of Theorem \ref{main theorem} combines ideas of
Gessel-Krattenthaler \cite[Proposition 1]{GKcyl} and Gasharov \cite[Theorem 3]{Gasharov}.

We begin by generalizing Young diagrams to non-partition shapes. For a composition $\alpha = (\alpha_1, \alpha_2,\dots,\alpha_l)$ the Young diagram of shape $\col(\alpha)$ is the set $\{(i,j): 1 \le j \le l, 1 \le i \le \alpha_j\}$ using matrix style coordinates. For a partition $\beta$ with $\col(\beta) \subset \col(\alpha)$, write $\col(\alpha, \beta)$ for $\col(\alpha) - \col(\beta)$. We define a $P$-array of shape $\col(\alpha, \beta)$ to be a function $A: \col(\alpha, \beta) \to P$ that is increasing in $P$ down columns. We write $P\text{-Array}(\alpha, \beta)$ for the set of $P$-arrays of shape $\col(\alpha, \beta)$.

As with $P$-tableaux, for a $P$-array $A$ let $A^d$ be the $P$-array obtained by gluing the first column of $A$ to the right of the last column of $A$ so that the copied column is shifted up by $d$ cells, i.e. $A^d$ has shape $\col(\bar{\alpha}, \bar{\beta})$ where $\bar{\alpha} = (\alpha_1 + d, \alpha_2 + d,\dots, \alpha_l + d, \alpha_1)$ and $\bar{\beta} = (\beta_1+d,\beta_2 + d,\dots,\beta_l+d,\beta_1)$.

\begin{exmp}
For $P$ the total order on $\bZ_{> 0}$,
\begin{equation}
\raisebox{-20pt}{$A =$ } \begin{ytableau}
\none & \none & \none & 1 \\
\none & 2 & 1 & 2 \\
3 & 3 & \none & 3 \\
5 & 4 \\
\none & 6 \\
\none & 7
\end{ytableau}
\end{equation}
is a $P$-array of shape $\col((4,6,2,3),(2,1,1))$ , and
\begin{eqnarray}
\raisebox{-20pt}{$A^2 =$ } \begin{ytableau}
\none & \none & \none & 1 & 3\\
\none & 2 & 1 & 2 & 5\\
3 & 3 & \none & 3 \\
5 & 4 \\
\none & 6 \\
\none & 7
\end{ytableau}
& \raisebox{-20pt}{ and } &
\raisebox{-20pt}{$A^3 =$ }
\begin{ytableau}
\none & \none & \none & \none & 3 \\
\none & \none & \none & 1 & 5\\
\none & 2 & 1 & 2 \\
3 & 3 & \none & 3 \\
5 & 4 \\
\none & 6 \\
\none & 7
\end{ytableau}.
\end{eqnarray}
\end{exmp}

Now fix a cylindric shape $\lambda/\mu/d$.
Write $\tilde{\bZ}^{\lambda_1}$ for the set of $\bk = (k_1,\dots, k_{\lambda_1}) \in \bZ^{\lambda_1}$ such that $k_1+k_2+\cdots+k_{\lambda_1} = 0.$
For a permutation $\pi \in S_{\lambda_1}$ and  $\bk \in \tilde{\bZ}^{\lambda_1}$, let $\pi_d^{\mathbf{k}}(\lambda)$ be the sequence of $\lambda_1$ integers given by
\begin{equation}
\pi_d^{\mathbf{k}}(\lambda)_i = k_i(\lambda_1+d) + \lambda_{\pi(i)}' + i - \pi(i)
\end{equation}
for $i \in [\lambda_1]$.
Note that when $\pi$ is the identity permutation and $\bk = \mathbf{0}$, $\pi_d^\bk(\lambda) = \lambda'$.

Let
\begin{equation}
\tilde{B} = \Big\{ (\pi, \bk, A) : \pi \in S_{\lambda_1}, \, \bk \in \tilde{\bZ}^{\lambda_1}, \,
A \in  P\text{-Array}(\col(\pi_d^\bk(\lambda), \mu'))  \Big\},
\end{equation}
\begin{equation}
B = \big\{(\pi, \bk , A) \in \tilde{B}: A^d \text{ is not a  $P$-tableau} \big\}.
\end{equation}
Note that if $\pi_d^\bk(\lambda)_i < \mu_i'$ for some $i$ and $\bk$, there are no $P$-arrays of shape $\pi_d^\bk(\lambda)$. Say that the sign of  $(\pi, \bk, A) \in \tilde{B}$ is the sign of  $\pi$.
We will give a sign-reversing involution $\Phi$ on $B$, therefore showing that
\begin{equation}
\label{e sum B}
\sum_{(\pi, \bk, A) \in B} \sgn(\pi)\mathbf{u}^A = 0,
\end{equation}
where  $\bu^A$ is the product of $u_a$'s over the entries  $a$ of  $A$.

\begin{defin}
For a $P$-array $A = [a_{s,t}]_{(s,t) \in \col(\alpha,\beta)}$ of shape $\col(\alpha, \beta)$, we say $a_{s,t}$ is \emph{empty} if $(s,t) \not\in \col(\alpha,\beta)$ and $s > \beta_t$. For columns $i,j$, $i<j$, we say column $i$ and column $j$ \emph{intersect} if there is some index $m$ so that $a_{m,j}$ is nonempty, $m-j+i+1 > \beta_i,$ and
\begin{enumerate}
\item $a_{m-j+i+1,i}$ is empty, or
\item $a_{m-j+i+1,i} >_P a_{m,j}$.
\end{enumerate}
We say $a_{m,j}$ is an \emph{intersection point} of columns $i$ and $j$, and that $a_{m-j+i+1,i}$ is a \emph{witness} to the intersection point.
\end{defin}

The following lemma applied to  $A^d$ shows that $B$ consists of all
triples  $(\pi, \bk, A) \in \tilde{B}$ such that $A^d$ has an intersection point.

\begin{lemma}
\label{lem tableau no intersect}
Let  $P$ be a \threeone-free poset and $A$ be  a $P$-array. The following are equivalent:
\begin{enumerate}
\item[(i)] $A$ is a $P$-tableau
\item[(ii)] $A$ has no adjacent intersecting columns
\item[(iii)] $A$ has no intersecting columns
\end{enumerate}

\end{lemma}
\begin{proof}

We first show that the definition of a $P$-tableau is the same as a $P$-array with no adjacent intersecting columns. A $P$-tableau is defined to be a $P$-array of skew partition shape that is nondecreasing in $P$ across rows. Now a $P$-array $A$ is of skew partition shape $\lambda/\mu$ only if there is no $c>1$ and $r>\mu_{i-1}$ so that $a_{r,c}$ is nonempty and $a_{r,c-1}$ is empty. But that is exactly condition (1) of the defintition of intersecting columns applied to $i=c-1$ and $j=c$. Likewise, the nondecreasing row condition of $P$-tableau is equivalent to saying that no adjacent columns satisfy condition (2) of the definition of intersecting columns. Therefore (i) and (ii) are equivalent.

Now if $A$ has adjacent intersecting columns, then $A$ has intersecting columns. So (iii) $\implies$ (ii).

It remains to show that if some columns $i$ and $j$ intersect, then there is some adjacent columns that intersect.
To do this, choose a pair $i,j$, $i < j$, of intersecting columns so that $j-i$ is as small as possible. We then show that $j-i = 1$. Suppose $j-i > 1$. Let $a_{m,j}$ be the intersection point between columns $i$ and $j$. Now $\beta_{j-1} \ge \beta_j$, so if $a_{m-1,j-1}$ is empty, columns $j-1$ and $j$ intersect and we are done. Therefore we assume that $a_{m-1,j-1}$ is nonempty. Now if $a_{m-j+i+1,i}$ is empty, then, as $a_{m-1,j-1}$ is nonempty, columns $i$ and $j-1$ intersect, which is a contradiction. We can then assume that $a_{m-j+i+1,i}$ is nonempty. As $m-j+i+1 \ge \beta_i \ge \beta_j$, we have that $a_{m-1,j}$ is nonempty, and $a_{m,j}$ is an intersection point between columns $i$ and $j$, we have $a_{m-1,j} <_P a_{m,j} <_P a_{m-j+i+1,i}$. If $a_{m-1,j} <_P a_{m-1,j-1}$ then columns $j-1$ and $j$ intersect, and if $a_{m-1,j-1} <_P a_{m-j+i+1,i}$ then column $i$ and $j-1$ intersect, which is a contradiction. However, if $a_{m-1,j} \not<_P a_{m-1,j-1}$ and $a_{m-1,j-1} \not<_P a_{m-j+i+1,i}$, then $a_{m-1,j}, a_{m,j}, a_{m-j+i+1,i}$, and $a_{m-1,j-1}$ form an induced \threeone, which contradicts the fact that $P$ has no induced $\threeone$ subposet.

Therefore, if $A$ is a $P$-array with a pair of intersecting columns, then $A$ must have an adjacent pair of intersecting columns, as desired.
\end{proof}

The next lemma shows that  $\tilde{B} \setminus B$ consists of exactly the triples
$(\text{id}, \mathbf{0}, A)$ with  $A$ a
cylindric $P$-tableau of shape $\lambda/\mu/d$.

\begin{lemma}
\label{lem:fixed points}
If $A$ is a $P$-array of shape $\col(\pi_d^\bk(\lambda), \mu')$ such that $A^d$ is a $P$-tableau, then $\pi$ is the identity permutation and $\bk = \mathbf{0}$.
\end{lemma}
\begin{proof}

Let $A$ be a $P$-array of shape $\col(\pi_d^\bk(\lambda), \mu')$ such that $A^d$ is a $P$-tableau.
By Lemma \ref{lem tableau no intersect}, for any $1 < i \le \lambda_1$, column $i$ does not intersect with column 1 or column $\lambda_1+1$. Therefore, from the first condition in the definition of intersecting columns, we have the following inequalities
\begin{equation}
\pi_d^\bk(\lambda)_1 \ge \pi_d^\bk(\lambda)_i - i + 2
\end{equation}
and
\begin{equation}
\pi_d^\bk(\lambda)_i + d  \ge \pi_d^\bk(\lambda)_1 - (\lambda_1 + 1) + i + 1
\end{equation}
which are equivalent to
\begin{equation}
\lambda_{\pi(1)}' - \lambda_{\pi(i)}' + \pi(i) - (\pi(1) + 1) \ge (k_i - k_1)(\lambda_1 + d)
\end{equation}
and
\begin{equation}
(k_i-k_1+1)(\lambda_1 + d) \ge \lambda_{\pi(1)}' - \lambda_{\pi(i)}' + \pi(i) - \pi(1) + 1.
\end{equation}

Now consider the case when $\pi(i) > \pi(1)$. Then we have that $0 \le \pi(i) - (\pi(1) + 1) < \lambda_1$ and $0 \le \lambda_{\pi(1)}' - \lambda_{\pi(i)}' \le \lambda_{1}' - \lambda_{\lambda_1}' \le d$. Therefore we have
\begin{equation}
\lambda_1 + d > \lambda_{\pi(1)}' - \lambda_{\pi(i)}' + \pi(i) - (\pi(1) + 1) \ge (k_i-k_1)(\lambda_1 + d),
\end{equation}
so $k_i - k_1 < 1$. Furthermore, we have
\begin{equation}
(k_i - k_1 + 1)(\lambda_1 + d) \ge \lambda_{\pi(1)}' - \lambda_{\pi(i)}' + \pi(i) - \pi(1) + 1 \ge 2,
\end{equation}
so $k_i - k_1 + 1 > 0$. Therefore we have that $k_i = k_1$.

Now consider the case when $\pi(i) < \pi(1)$. Then $0 \ge \pi(i) - (\pi(1) + 1) \ge -\lambda_1$ and $0 \ge \lambda_{\pi(1)}' - \lambda_{\pi(i)}' \ge \lambda_{\lambda_1}' - \lambda_{1}' \ge -d.$
Then
\begin{equation}
0 \ge \lambda_{\pi(1)}' - \lambda_{\pi(i)}' + \pi(i) - (\pi(1) + 1) \ge (k_i-k_1)(\lambda_1 + d),
\end{equation}
so $k_i < k_1$, and
\begin{equation}
(k_i - k_1 + 1)(\lambda_1 + d) \ge \lambda_{\pi(1)}' - \lambda_{\pi(i)}' + \pi(i) - \pi(1) + 1 \ge -\lambda_1 - d + 1,
\end{equation}
so $k_i - k_1 + 1 > -1.$ Therefore we have that $k_1 = k_i + 1$.

Then
\begin{equation}
k_1 + k_2 + \cdots + k_{\lambda_1} = \lambda_1 k_1 - \#\{i : \pi(i) < \pi(1)\} = 0.
\end{equation}
But $0 \le \#\{i : \pi(i) < \pi(1)\} < \lambda_1$, so $k_1 = k_2 = \cdots = k_{\lambda_1} = 0.$
Now we show that $\pi$ must be the identity permutation. With $\bk = \mathbf{0}$, we have that if $i < j$ and $\pi(i) > \pi(j)$ then $\lambda_{\pi(i)}' \le \lambda_{\pi(j)}'$ and
\begin{equation}
\pi_d^\bk(\lambda)_i = \lambda_{\pi(i)}' + i - \pi(i) < \lambda_{\pi(j)}' + j - \pi(j) -j +i +1 = \pi_d^\bk(\lambda)_j - j + i + 1,
\end{equation}
so columns $i$ and $j$ intersect. Therefore we have that $\pi(i) < \pi(j)$ whenever $i < j$, so $\pi$ must be the identity permutation.
\end{proof}

\begin{defin}
Let $\alpha = (\alpha_1,\alpha_2,...,\alpha_l)$ be a composition and $\beta \subset \alpha$ be a partition.
For a $P$-array $A$ of shape $\col(\alpha, \beta)$ with columns $i$ and $j$ intersecting, we define a swap at $i$ and $j$ by letting $\swap(A,i,j)$ be the $P$-array obtained in the following way. Let $a_{m,j}$ be the $P$-minimal intersection point between columns $i$ and $j$, and let $a_{m-j+i+1,i}$ be the witness to this intersection.
Let
\begin{eqnarray}
R_{i,j} = \{a_{r,j}: r > m\} & \text{and} &  L_{i,j} = \{a_{r,i}: r \ge m-j+i+1\}.
\end{eqnarray}
Then $\swap(A, i, j)$ is the $P$-array obtained by moving $R_{i,j}$ to column $i$, $L_{i,j}$ to column $j$, and fixing the rest of the cells.

Given an integer $d$, $\beta_1 - \beta_l \le d$, 
and columns $i$ and $j$ so that $i$ and $j$ intersect in $A^d$, we define the cylindric swap at $i$ and $j$ by letting $\swap(A,i,j,d)$ be the $P$-array obtained in the following way. If $j \neq l + 1$, $\swap(A,i,j,d) = \swap(A,i,j)$. 
Now if $j = l + 1$, we again let $a_{m,j}$ be the $P$-minimal intersection point between columns $i$ and $j$, and take
\begin{eqnarray}
R_{i,j} = \{a_{r,j}: r > m\} = \{a_{r,1}: r > m+d\} & \text{and} & L_{i,j} = \{a_{r,i}: r \ge m-j+i+1\}.
\end{eqnarray}
Then $\swap(A,i,j,d)$ is the $P$-array obtained from  $A$ by moving $R_{i,j}$ to column $i$, $L_{i,j}$ to column 1, and fixing the rest of the cells of $A$.

Note that the cylindric swap is an operation on $A$ rather than on $A^d$, even though we use an intersection point of $A^d$ to define it.
\end{defin}

\begin{exmp}
Let $P$ be the poset on $\bZ_{>0}$ so that $i <_P j$ if $i<j$ and $j-i > 1$, and consider the following $P$-arrays $A$ and  $A^1$
\begin{equation}
A= \raisebox{20pt}{$\begin{ytableau}
\none & 2 & 3   \\
1 & 4 & 5   \\
4 & 9 & 7 \\
6
\end{ytableau}$}\, ,
\quad \quad \qquad
A^1= \raisebox{20pt}{ $\begin{ytableau}
\none & 2 & 3 & \underline{1} \\
1 & 4 & 5 & 4 \\
4 & 9 & \underline{7} & 6 \\
6
\end{ytableau}$}\, .
\end{equation}
The underlined entries indicate the intersection points of  $A^1$.

We then have
\begin{eqnarray}
\swap(A,2,3,1) = \raisebox{20pt}{$ \begin{ytableau}
\none & 2 & 3 \\
1 & 4 & 5 \\
4 & \none & 7 \\
6 & \none & 9
\end{ytableau}$}
& \ \ \ \text{and} \ \ \  &
\swap(A, 3,4,1) = \raisebox{20pt}{$\begin{ytableau}
\none & 2 & 4 \\
1 & 4 & 6 \\
3 & 9 \\
5 \\
7
\end{ytableau}$}\, .
\end{eqnarray}
\end{exmp}

\begin{lemma}
\label{lem: swap}
If $A$ is a $P$-array of shape $\col(\pi_d^{\mathbf{k}}(\lambda), \mu')$ so that columns $i$ and $j$ intersect in $A^d$, then $\swap(A, i, j, d)$ is a $P$-array of shape $\col(\sigma_{d}^{\mathbf{k}'}(\lambda), \mu')$ where
\begin{equation}
\sigma = \begin{cases} \pi \circ (i \ j) & \text{if } j \neq \lambda_1 + 1 \\
\pi \circ (1 \ i) & \text{if } j = \lambda_1 + 1
\end{cases}
\end{equation}
and
\begin{equation}
\bk' = \begin{cases}
(k_1,\dots,k_j,\dots,k_i,\dots,k_{\lambda_1}) & \text{if } j \neq \lambda_1 + 1 \\
(k_i + 1, k_2,\dots,k_{i-1},k_1-1,k_{i+1},\dots,k_{\lambda_1}) & \text{if } j = \lambda_1 + 1
\end{cases}.
\end{equation}
\end{lemma}

\begin{proof}
Let $A$ be a $P$-array of shape $\col(\pi_d^k(\lambda), \mu')$ such that columns $i$ and $j$ of $A^d$ intersect. We first consider the case when $1 \le i < j \le \lambda_1$. We then have that $\swap(A, i, j)$  has shape $\col(\alpha, \mu')$ where $\alpha_i = \pi_d^\bk(\lambda)_j +i -j$, $\alpha_j = \pi_d^\bk(\alpha)_i + j -i$, and the rest of the columns are fixed. So we have
\begin{equation}
\alpha_i = \lambda_{\pi(j)}' + j - \pi(j) + k_j(\lambda_1 + d) + i - j = \lambda_{\pi(j)}' + i - \pi(j) + k_j(\lambda_1 + d),
\end{equation}
and
\begin{equation}
\alpha_j = \lambda_{\pi(i)}' + i - \pi(i) + k_i(\lambda_1 + d) + j - i = \lambda_{\pi(i)}' + j - \pi(i) + k_i(\lambda_1 + d).
\end{equation}
So we have $\alpha = \sigma_d^{\bk'}(\lambda)$ where $\sigma = \pi \circ (i\ j)$ and $\bk' = (k_1,\dots,k_j,\dots,k_i,\dots,k_{\lambda_1})$ as desired.

Now consider the case when $j = \lambda_1 + 1$. Then $\swap(A,i,j,d)$ has shape $\col(\alpha, \mu')$ where
\begin{align}
\alpha_i &= \pi_d^\bk(\lambda)_1 - (\lambda_1 + 1) + i- d = \lambda_{\pi(1)}' + 1 - \pi(1) + k_1(\lambda_1 + d) - (\lambda_1 + 1) + i - d \\
&= \lambda_{\pi(1)}' + i - \pi(1) + (k_1 - 1)(\lambda_1 + d),
\end{align}
\begin{align}
\alpha_1 &= \pi_d^\bk(\lambda)_i - i + (\lambda_1 + 1) + d = \lambda_{\pi(i)}' + i - \pi(i) + k_i(\lambda_1 + d) - i + (\lambda_1 + 1) + d\\
 &= \lambda_{\pi(i)}' + 1 - \pi(i) + (k_i + 1)(\lambda_1 + d),
\end{align}
and the rest of the columns are fixed.
Therefore we have that $\alpha = \sigma_d^{\bk'}(\lambda)$ where $\sigma = \pi \circ (1\ i)$ and $\bk' = (k_i + 1, k_2,\dots, k_{i-1}, k_1 - 1, k_{i+1},\dots,k_{\lambda_1})$ as desired.
\end{proof}

We define the involution $\Phi$ on  $B$ in the following way:
for a $P$-array $A$, let $a_{m,j}$ be the rightmost $P$-minimal intersection point in $A^d$.  Let $i$ be the rightmost column index so that columns $i$ and $j$ intersect with intersection point $a_{m,j}$. We let
 $\phi(A) = \swap(A, i, j, d)$ and
$\Phi((\pi, \bk, A)) = (\sigma, \bk', \phi(A))$, where  $\sigma$ and  $\bk'$ are as in
Lemma \ref{lem: swap}.

Lemma \ref{lem: swap} shows that  $\Phi$ is  sign-reversing.
To show that  $\Phi$ is an involution, we need to show that the rightmost $P$-minimal intersection of $A^d$ is also the rightmost $P$-minimal intersection point of $\phi(A)^d$.
This is established in the next two lemmas.

\begin{lemma}
\label{lem: inv1}
If $a_{m,j}$ is the rightmost $P$-minimal intersection point of $A^d$, then $a_{m,j}$ is a $P$-minimal intersection point of $\phi(A)^d$. Furthermore, if $i$ is the rightmost column so that $a_{m,j}$ is an intersection point between columns $i$ and $j$, $i$ is the rightmost column in $\phi(A)^d$ so that $a_{m,j}$ is an intersection point between columns $i$ and $j$.
\end{lemma}
\begin{proof}
First we show that $a_{m,j}$ is an intersection point in $\phi(A)^d$. As $\phi$ swaps the positions of $\{a_{m-j+i+1,i},a_{m-j+i+2,i},\dots,a_{\alpha_i, i}\}$ and $\{a_{m+1,j}, a_{m+2,j},\dots,a_{\alpha_j,j}\}$, and either $a_{m,j} <_P a_{m+1,j}$ or $a_{m+1,j}$ is not defined, $a_{m,j}$ remains a intersection point of columns $i$ and $j$. As $\phi$ fixes the entries above row $m$ of column $j$ and the entries above $m-j+i$ in column $i$, $a_{m,j}$ remains the $P$-minimal intersection point between columns $i$ and $j$. As $\phi$ fixes columns $\{i+1,i+2,\dots,j-1\}$, column $i$ is the rightmost column so that $a_{m,j}$ is an intersection point.

Now we will show that $a_{m,j}$ is $P$-minimal among intersection points of $\phi(A)^d$. Consider some intersection point $x$ of $\phi(A)^d$ such that $x <_P a_{m,j}.$ As each entry that is different in $\phi(A)^d$ from the corresponding entry of $A$ is greater than $a_{m,j}$ in $P$, each such entry is also greater than $x$ in $P$. Therefore, $x$ relates to each entry of $A^d$ in the same way as it relates to the corresponding entry of $\phi(A)^d$. So, if $x$ is an intersection point of $\phi(A)^d$ it must have also been an intersection point of $A^d$, which contradicts the fact that $a_{m,j}$ is a $P$-minimal intersection point. Therefore $a_{m,j}$ is a $P$-minimal intersection point in both $A^d$ and $\phi(A)^d$ as desired.
\end{proof}
\begin{lemma}
\label{lem: inv2}
If $a_{m,j}$ is the rightmost $P$-minimal intersection point of $A^d$, $a_{m,j}$ is the rightmost $P$-minimal intersection point of $\phi(A)^d$.
\end{lemma}
\begin{proof}
Suppose $x \neq a_{m,j}$ is the rightmost $P$-minimal intersection point of $\phi(A)^d$. Then $x$ is in some column $k$ with $j < k \le \lambda_1+1$.
As every entry that is moved by $\phi$ is greater than $a_{m,j}$ in $P$, $x$ must be in the same location in $A^d$ and $\phi(A)^d$. Next observe that if $y$ is a witness to $x$ being an intersection point, $y$ must be an entry of the southeast diagonal directly below the diagonal containing $x$. e.g. the bullets in the following array indicate the possible witnesses of $x$ being an intersection point.
\begin{equation}
\begin{ytableau}
\bullet & & & &\\
& \bullet & & & \\
& & \none[\bullet] & x & \\
& \none & \none &
\end{ytableau}
\end{equation}
Now as $j \neq \lambda_1+1$, $\phi$ fixes southeast diagonals in columns $1,2,\dots,\lambda_1$. Therefore if $x$ is an intersection point to the right of column $j$ and incomparable to $a_{m,j}$ in $\phi(A)^d$, $x$ must have been an intersection point in $A^d$. But this is a contradiction.
\end{proof}

We now combine the previous lemmas to prove Theorem \ref{main theorem}.
\begin{proof}
We wish to show that
\begin{equation}
s_{\lambda/\mu/d}^P(\mathbf{u}) = \sum_{\bk \in \tilde{\bZ}^{\lambda_1}}
\det [ e_{k_i(\lambda_1 + d) + \lambda_i' - \mu_j' - i + j}^P(\mathbf{u})] = \sum_{T \in \CT_P(\lambda/\mu/d)} \mathbf{u}^T.
\end{equation}
By Lemmas \ref{lem: swap}, \ref{lem: inv1}, and \ref{lem: inv2},  $\Phi$ is a sign-reversing involution,
and by Lemma \ref{lem:fixed points},
$\tilde{B} \setminus B$ consists of exactly the triples
$(\text{id}, \mathbf{0}, A)$ with  $A$ a
cylindric $P$-tableau of shape $\lambda/\mu/d$.
Hence
\begin{equation}
\sum_{T \in \CT_P(\lambda/\mu/d)} \mathbf{u}^T = \sum_{\bk \in \tilde{\bZ}^{\lambda_1}} \sum_{\pi \in S_{\lambda_1}} \mathrm{sgn}(\pi) \sum_{T \in P\text{-Array}(\pi_d^\bk(\lambda), \mu')} \mathbf{u}^T.
\end{equation}
Now for a composition $\alpha$ and partition $\beta \subset \alpha$, we have
\begin{equation}
\sum_{T \in P\text{-Array}(\alpha, \beta)}\mathbf{u}^T = \prod_{i = 1}^{\ell(\alpha)} e_{\alpha_i - \beta_i}^P(\mathbf{u}),
\end{equation}
so
\begin{equation}
\sum_{T \in \CT_P(\lambda/\mu/d)} \mathbf{u}^T = \sum_{\bk \in \tilde{\bZ}^{\lambda_1}} \sum_{\pi \in S_{\lambda_1}} \mathrm{sgn}(\pi) e_{k_i(\lambda_1+d) + \lambda_{\pi(i)}' + i - \pi(i) - \mu_i'}^P(\mathbf{u}).
\end{equation}
Permutation of the indices fixes $\tilde{\bZ}^{\lambda_1}$, so we have
\begin{align*}
\sum_{\bk \in \tilde{\bZ}^{\lambda_1}} \sum_{\pi \in S_{\lambda_1}} \mathrm{sgn}(\pi)& e_{k_i(\lambda_1+d) + \lambda_{\pi(i)}' + i - \pi(i) - \mu_i'}^P(\mathbf{u}) = \sum_{\pi \in S_{\lambda_1}} \sum_{\bk \in \tilde{\bZ}^{\lambda_1}} \mathrm{sgn}(\pi) e_{k_{\pi(i)}(\lambda_1+d) + \lambda_{\pi(i)}' + i - \pi(i) - \mu_i'}^P(\mathbf{u})\\
&= \sum_{\bk \in \tilde{\bZ}^{\lambda_1}} \sum_{\pi \in S_{\lambda_1}} \mathrm{sgn}(\pi)e_{k_{\pi(i)}(\lambda_1+d) + \lambda_{\pi(i)}' + i - \pi(i) - \mu_i'}^P(\mathbf{u}).
\end{align*}
As the determinant of an $n \times n$ matrix $X$ is given by $\det X = \sum_{\pi \in S_n} \sgn(\pi) \prod_{i=1}^n X_{\pi(i),i}$, we have
\begin{equation}
\sum_{T \in \CT_P(\lambda/\mu/d)} \mathbf{u}^T = \sum_{\bk \in \tilde{\bZ}^{\lambda_1}} \det [ e_{k_i(\lambda_1 + d)+ \lambda_i' - \mu_j' -i + j }^P(\mathbf{u})],
\end{equation}
as desired.
\end{proof}

\begin{remark}
In the proof of Theorem \ref{main theorem}, we used the fact that $S_{\lambda_1}$ acts on $\tilde{\bZ}^{\lambda_1}$ to show that
\begin{equation*}
\sum_{\bk \in \tilde{\bZ}^{\lambda_1}} \sum_{\pi \in S_{\lambda_1}} \mathrm{sgn}(\pi) e_{k_i(\lambda_1+d) + \lambda_{\pi(i)}' + i - \pi(i) - \mu_i'}^P(\mathbf{u}) = \sum_{\bk \in \tilde{\bZ}^{\lambda_1}} \sum_{\pi \in S_{\lambda_1}} \mathrm{sgn}(\pi)e_{k_{\pi(i)}(\lambda_1+d) + \lambda_{\pi(i)}' + i - \pi(i) - \mu_i'}^P(\mathbf{u}).
\end{equation*}
Therefore we can obtain two equivalent determinantal formulas for $s_{\lambda/\mu/d}^P(\mathbf{u})$ by replacing $k_i$ with $k_j$:
\begin{equation}
s_{\lambda/\mu/d}^P(\mathbf{u}) = \sum_{\bk \in \tilde{\bZ}^{\lambda_1}} \det [e_{k_i(\lambda_1 + d) + \lambda_i' - \mu_j' -i + j }^P(\mathbf{u})] = \sum_{\bk \in \tilde{\bZ}^{\lambda_1}} \det [ e_{k_j(\lambda_1 + d)+ \lambda_i' - \mu_j' -i + j  }^P(\mathbf{u})].
\end{equation}
\end{remark}

\section{Special Cases}
\label{s special cases}

We now examine two special cases of cylindric Schur functions with particularly nice monomial expansions.

Take the partition $\lambda = (r^c)$ whose Young diagram is a rectangle with width $r$ and height $c$. Then for $T \in \CSSYT(\lambda/ \varnothing/ 0),$ if $a_1, a_2,\dots, a_r$ are the entries of a row in $T$, then
\begin{equation}
a_1 \le a_2 \le \cdots \le a_r \le a_1,
\end{equation}
so $a_1 = a_2 = \cdots = a_r$.
Therefore we have
\begin{equation}
\label{e cylind schur 0}
s_{(r^c)/\varnothing/0}(\mathbf{x}) = m_{(r^c)}(\mathbf{x}).
\end{equation}

If $\lambda = (r^ct) \vdash n$ is the partition with $c$ rows of length $r$ and another row of length $t \le r$, then any cylindric tableau  $T$ of shape $\lambda/\varnothing/1$ has the following form:
 the word  $a_1 \cdots a_n$ formed by concatenating the rows of  $T$ from top to bottom must satisfy
 $a_1 \le a_2 \le\cdots\le a_n$ and  $a_i < a_{i+r}$ for all  $1 \le i \le n-r$.
Therefore, given any multiset of positive integers so that no number is repeated more than $r$ times, we can fill the diagram of shape $\lambda/\varnothing/1$ in exactly one way. Hence
\begin{equation}
\label{e cylind schur 1}
s_{(r^c t)/ \varnothing/ 1}(\mathbf{x}) = \sum_{\lambda_1 \le r} m_{\lambda}(\mathbf{x}).
\end{equation}

The monomial expansions \eqref{e cylind schur 0} and \eqref{e cylind schur 1} combined with Corollary \ref{c sum of e coefs} yield the following results.

\begin{corollary}
\label{c: rect pos}
For a \threeone-free poset $P$, the coefficient of $e_{(r^c)}(\mathbf{x})$ in the $e$-expansion of $X_{\inc(P)}(\mathbf{x})$ is the number of standard cylindric $P$-tableaux of shape $(r^c)/\varnothing/0$.
\end{corollary}
\begin{corollary}
\label{c: part sum}
For a \threeone-free poset $P$ and positive integer  $r$,
letting $c^P_\lambda$ denote the coefficients
in  $X_{\inc(P)}(\mathbf{x}) = \sum_\lambda c^P_\lambda e_\lambda(\mathbf{x})$, we have
\begin{equation}
\sum_{\substack{\lambda \vdash |P| \\ \lambda_1 \le r}} c^P_{\lambda} =
\# \big\{T \in \CT_P((r^c t)/\varnothing/1) : \text{ $T$ is standard}\big\},
\end{equation}
where $(r^c t)$ is the partition of $|P|$ with the maximal number of rows of length $r$.
\end{corollary}

Corollary \ref{c: rect pos} recovers a result of Clearman-Hyatt-Shelton-Skandera \cite[Theorem 4.7 (v-b)]{CHSS} which gave a combinatorial interpretation to a theorem of Stembridge \cite[Theorem 2.8]{StembridgeImmanant}.

Corollary \ref{c: part sum} is reminiscent of the following theorem of Stanley \cite{Stanleychromatic}:
\begin{theorem}
Let $c^G_\lambda$ be the coefficient of $e_\lambda$ in a chromatic symmetric function $X_G$. Then
\begin{equation}
\sum_{\lambda\, :\, \ell(\lambda) = j} c^G_\lambda = \text{\rm sink}(G, j)
\end{equation}
where $\text{\rm sink}(G,j)$ is the number of acyclic orientations of $G$ with $j$ sinks.
\end{theorem}

\section{Acknowledgements}
The author would like to thank Jonah Blasiak for his mentorship and Holden Eriksson for helpful discussions. This project was started as part of the Co-op program at Drexel University.
The author was supported by NSF Grant DMS-1855784.

\bibliographystyle{plain}
\bibliography{mycitations}

\end{document}